\documentclass{amsart}[11pt]
\usepackage{amsfonts,mathrsfs,amsmath,amssymb}
\numberwithin{equation}{section}
\newtheorem{theo}[equation]{Theorem}
\newtheorem{pro}[equation]{Proposition}
\newtheorem{cor}[equation]{Corollary}
\newtheorem{lem}[equation]{Lemma}
\newtheorem{defi}[equation]{Definition}
\newtheorem{rem}[equation]{Remark}

\def\vn{\varepsilon}
\def\ot{\otimes}
\def\om{\omega}
\def\mg{\mathfrak{g}}

\def\lan{\langle}
\def\ran{\rangle}
\def\al{\alpha}

\def\be{\beta}
\def\De{\Delta}
\def\ga{\gamma}

\def\ep{\epsilon}
\def\de{\delta}
\def\pa{\partial}

\def\D{{\mathcal D}}

\def\bZ{{\mathbb Z}}
\def\bQ{{\mathbb Q}}

\def\bC{{\mathbb C}}
\def\lra{\longrightarrow}

\begin{document}

\title[Notes on two-parameter quantum groups, (I)]
{Notes on two-parameter quantum groups, (I) }

\author[Hu]{Naihong Hu$^\star$}
\address{Department of Mathematics, East China Normal University,
Shanghai 200062, PR China} \email{nhhu@euler.math.ecnu.edu.cn}

\author[Pei]{Yufeng Pei}
\address{Department of Mathematics, East China Normal University,
Shanghai 200062, PR China}\email{peiyufeng@gmail.com}
\thanks{$^\star$N.H., Corresponding author, supported by the NNSF (Grant 10431040),
the TRAPOYT, the FUDP and the Priority Academic Discipline from the
MOE of China, the SRSTP from the STCSM, the Shanghai Priority
Academic Discipline from the SMEC}


\begin{abstract}
A simpler definition for a class of two-parameter quantum groups
associated to semisimple Lie algebras is given in terms of Euler
form. Their positive parts turn out to be $2$-cocycle deformations
of each other under some conditions. An operator realization of the
positive part is given.
\end{abstract}
\maketitle

\section{Introduction}
The notion of quantum groups was introduced by V. Drinfel'd and M.
Jimbo, independently, around 1985 in their study of the quantum
Yang-Baxter equations. Quantum groups $U_q(\mg)$, depending on a
single parameter $q$, are certain families of Hopf algebras that are
deformations of universal enveloping algebras of symmetrizable
Kac-Moody algebras. In the early 90¡¯s of the last century, much
work had been done on their multiparameter generalizations, which
can be obtained  by twisting the algebra structure via a 2-cocycle
on an indexed free abelian group (see \cite{AST}) or by twisting the
coalgebra structure in the spirit of Drinfeld (see \cite{Res},
\cite{HW}). Note that a $2$-cocycle (or a Drinfeld twist)
deformation is an important method to yield new (twisted) bialgebras
from old ones.

\smallskip
Motivated by the work on down-up algebras \cite{BW0}, Benkart and
Witherspoon, et al \cite{BW1, BW2, BW3, BKL} investigated the
two-parameter quantum groups of the general linear Lie algebra
$\frak{gl}_n$ and the special linear Lie
algebra $\frak{sl}_n$. 
Later on, Bergeron, Gao and Hu \cite{BGH1, BGH2} developed the
corresponding theory for two-parameter quantum orthogonal and
symplectic groups. Recently, Hu et al continued this project (see
\cite{HS, BH} for exceptional types $G, E$, \cite{HW2}, \cite{CHW},
\cite{BH2} for restricted types $B$, $C$, $D$, and \cite{HRZ},
\cite{HZ} for untwisted affine types, etc.).

\smallskip
In this note, we give a simpler definition for a class of
two-parameter quantum groups $U_{r,s}(\mg)$ associated to semisimple
Lie algebras in terms of the Euler form (or say, Ringel form). As in
\cite{BW1, BGH1, HS, BH, HRZ, HZ}, these quantum groups also possess
the Drinfel'd double structures and the triangular decompositions
(see Section 2). As a main point of this note, we show that the
positive parts of quantum groups under consideration are $2$-cocycle
deformations of each other as $Q^+$-graded associative
$\bC$-algebras if the parameters satisfy certain conditions (see
Section 3). This affords an insight into the interrelation between
the two-parameter quantum groups we defined and the one-parameter
Drinfeld-Jimbo ones. In Section 4, we get an operator realization of
the positive part of $U_{r,s}(\mg)$ by assigning the canonical
generators $e_i$'s with some skew differential operators in the
sense of Kashiwara (\cite{Ka}).

\bigskip
\section{Euler form and definition of two-parameter quantum groups}
\medskip

 Let $C=(a_{ij})_{i,j\in I}$ be a Cartan matrix of finite type and
 $\mg$ the associated semisimple Lie algebra over $\bQ$. Let $\{\,d_i\mid i\in I\,\}$ be
 a set of relatively prime
positive integers such that $d_ia_{ij}=d_{j}a_{ji}$ for $i,j\in I$.
Let $\bQ(r,s)$ be the function field in two variables $r,\,s$ over
the field $\bQ$ of rational numbers. Denote $r_i=r^{d_i},
\,s_i=s^{d_i}$ for $i\in I$.

Let $\lan-,-\ran$ be the bilinear form, which is called the Euler
form (or Ringel form), on the root lattice $Q$ defined by
$$
\lan i,j\ran:=\lan\al_i,\al_j\ran=
\begin{cases}
&d_ia_{ij}\quad i<j,\\
&d_i\quad\quad\, i=j,\\
&0\quad\quad\ \   i>j.
\end{cases}
$$

\begin{defi}
The two-parameter quantum group $U_{r,s}(\mg)$ is a unital
associative algebra over $\bQ(r,s)$ generated by
$e_i,f_i,\om_i^{\pm1},\om_i'^{\pm1},$ $i\in I$, subject to the
following relations$\,:$
\begin{eqnarray*}
 &(R1)&\quad
\om_i^{\pm1}\om_j^{\pm1}=\om_j^{\pm1}\om_i^{\pm1},
\qquad\qquad\ \om_i'^{\pm1}\om_j'^{\pm1}=\om_j'^{\pm1}\om_i'^{\pm1},\\
& &\quad\om_i^{\pm1}\om_j'^{\pm1}=\om_j'^{\pm1}\om_i^{\pm1},\qquad
\qquad\om_i^{\pm1}\om_i^{\mp1}=\om_i'^{\pm1}\om_i'^{\mp1}=1.\\
&(R2)&\quad \om_i\,e_j\,\om_i^{-1}=r^{\lan j, i\ran}s^{-\lan i,
j\ran} e_j, \qquad\qquad\om_i'\,e_j\,\om_i'^{-1}=r^{-\lan
i,j\ran}s^{\lan j,i\ran} e_j.
\\
&(R3)&\quad\om_i\,f_j\,\om_i^{-1}=r^{-\lan j, i\ran}s^{\lan i,
j\ran} f_j,\qquad\qquad
\om_i'\,f_j\,\om_i'^{-1}=r^{\lan i,j\ran}s^{-\lan j,i\ran}f_j.\\
&(R4)&\quad e_if_j-f_je_i=\delta_{i,j}\frac{\om_i-\om_i'}{r_i-s_i}.\\
&(R5)&\quad\sum_{k=0}^{1-a_{ij}} (-1)^{k}
\binom{1-a_{ij}}{k}_{r_is_i^{-1}} c^{(k)}_{ij}
e_{i}^{1-a_{ij}-k}\,e_{j}\,
e_{i}^{k} =0, \qquad\,\; (i\neq j),\\
&(R6)&\quad \sum_{k=0}^{1-a_{ij}} (-1)^{k}
\binom{1-a_{ij}}{k}_{r_is_i^{-1}} c^{(k)}_{ij}\, f_{i}^{k}\,f_{j}\,
f_{i}^{1-a_{ij}-k} =0, \qquad (i\neq j),
\end{eqnarray*}
where $c^{(k)}_{ij}=(r_is_i^{-1})^{\frac{k(k-1)}{2}} r^{k\lan
j,i\ran}s^{-k\lan i,j\ran},\ \textit{for } \ i\neq j$, and for a
symbol $v$, we set the notations$\,:$
\begin{gather*}
(n)_v=\frac{v^n-1}{v-1},\qquad (n)_v!=(1)_v(2)_v\cdots(n)_v,
\\
\binom{n}{k}_v=\frac{(n)_v!}{(k)_v!(n-k)_v!},\qquad\textit{for }
n\ge k\ge 0,
\end{gather*}
and $(0)_v!=1$.
\end{defi}
The algebra $U_{r,s}(\mg)$ has a Hopf algebra structure
with the comultiplication, the counit and the antipode given by:
\begin{gather*}
\Delta(\om_i^{\pm1})=\om_i^{\pm1}\ot\om_i^{\pm1}, \qquad
\Delta({\om_i'}^{\pm1})={\om_i'}^{\pm1}\ot{\om_i'}^{\pm1},\\
\Delta(e_i)=e_i\ot 1+\om_i\ot e_i, \qquad \Delta(f_i)=1\ot
f_i+f_i\ot \om_i',\\
\vn(\om_i^{\pm1})=\vn({\om_i'}^{\pm1})=1, \qquad\qquad
\vn(e_i)=\vn(f_i)=0,\\
S(\om_i^{\pm1})=\om_i^{\mp1}, \qquad\quad\qquad
S({\om_i'}^{\pm1})={\om_i'}^{\mp1},\\
S(e_i)=-\om_i^{-1}e_i,\qquad S(f_i)=-f_i\,{\om_i'}^{-1}.
\end{gather*}

\begin{rem}
$($\textrm{\rm i}$)$\  Let $r=q$, $s=q^{-1}$. Then $U_{q,q^{-1}}$
modulo the Hopf ideal generated by $\om_i'-\om_i^{-1}$ $(i\in I)$ is
isomorphic to the standard one-parameter quantum group $U_{q}(\mg)$
defined by Drinfel'd and Jimbo $($see \cite{Ja}$)$.

$($\textrm{\rm ii}$)$ \ Let $r=q^2$, $s=1$. Then $U_{q^2,1}^+$ is
isomorphic to the $($nontwisted$)$ generic Hall algebra introduced
by Ringel $($\cite{Ri}$)$.

$($\textrm{\rm iii}$)$ \ For the type $D$ case, the definition above
is distinct from that given in \cite{BGH1}.

$($\textrm{\rm iv}$)$ \ Definition 2.1 might be adopted to define
the affine cases, but the resulting quantum groups in this fashion
$($because of no Drinfeld realization to be found in these cases$)$
are different from that given in \cite{HRZ}.

\end{rem}

Let $U_{r,s}^+$ (resp., $U_{r,s}^-$) be the subalgebra of
$U_{r,s}:=U_{r,s}(\frak g)$ generated by the elements $e_i$ (resp.,
$f_i$) for $i\in I$, and $U^0$ the subalgebra of $U_{r,s}$ generated
by $\om_i^{\pm1},\om_i'^{\pm1}$  for $i\in I$. Moreover, let
$U_{r,s}^{\geq0}$ (resp., $U_{r,s}^{\leq0}$) be the subalgebra of
$U_{r,s}$ generated by the elements $e_i,\,\om_i^{\pm1}$ for $i\in
I$ (resp., $f_i,\, \om_i'^{\pm1}$ for $i\in I$). For each $\mu\in Q$
(the root lattice of $\frak g$), we define elements $\om_{\mu}$ and
$\om_{\mu}'$ by
$$
\om_{\mu}=\prod_{i\in I}\om_{i}^{\mu_i},\quad \om'_{\mu}=\prod_{i\in
I}{\om_{i}'}^{\mu_i},\qquad\textit{for } \ \mu=\sum_{i\in
I}\mu_i\al_i\in Q.
$$
 For $\be\in Q^+$ (a fixed positive root lattice), let
$$
U_{r,s}^{\pm\be}=\left\{x\in
U_{r,s}^{\pm}\,\left|\,\om_{\mu}x\om_{-\mu}
=r^{\lan\be,\mu\ran}s^{-\lan\mu,\be\ran}x,\,
\om_{\mu}'x\om_{-\mu}'
=r^{-\lan\mu,\be\ran}s^{\lan\be,\mu\ran}x,\, \forall\ \mu\in
Q\right\}\right.,
$$
then $U_{r,s}^{\pm}=\bigoplus_{\be\in Q^+}U_{r,s}^{\pm\be}$ are
$Q^+$-graded.

\begin{pro}
For any $i\in I$, we have the $\bQ$-algebra automorphism $\Phi$ and
the $\bQ(r,s)$-algebra anti-automorphism $($or say, involution$)$
$\Psi$ of $U_{r,s}(\mg)$ defined by
\begin{gather*}
\Phi(r)=s^{-1},\  \Phi(s)=r^{-1},\   \Phi(e_i)=f_i,\
\Phi(f_i)=r_is_ie_i,\ \Phi(\om_i)=\om_i',\ \Phi(\om_i')=\om_i.\\
\Psi(e_i)=f_i,\quad\Psi(f_i)=e_i,\quad\Psi(\om_i)=\om_i,\quad\Psi(\om_i')=\om_i'.
\end{gather*}
\end{pro}
\begin{proof}
 It is straightforward to check the statements.
\end{proof}

\begin{pro} There exists a unique  bilinear skew pairing
$$(\,,):\, U_{r,s}^{\leq0}\times U_{r,s}^{\geq0}\lra \bQ(r,s)$$
such that for all $x,\,x'\in U_{r,s}^{\geq0}$, \ $y,\,y'\in
U_{r,s}^{\leq0}$, \
 $\mu,\,\nu\in Q$, and \,$i,\,j\in I$,
\begin{equation*}
\begin{split}
( y,\,xx' )&=(\De(y),\,x'\ot x),\\
(yy',\,x) &=(y\ot y',\,\De(x)),\\
(f_i,\,e_j)&=\delta_{i,j}\frac{1}{s_i-r_i},\\
(\om_{\mu}',\,\om_{\nu})&=r^{\lan\mu,\nu\ran}s^{-\lan\nu,\mu\ran},\\
(\om_{\mu}',\,e_i)&=0,\\
(f_i,\,\om_{\mu}) &=0.
\end{split}
\end{equation*}
\end{pro}
\begin{proof}
The coalgebra structure of $U_{r,s}(\mg)$ defines an algebra
structure on $ (U_{r,s}^{\geq0})^* $ by
$$
(\gamma_1\ga_2)(x):=(\ga_1\ot\ga_2)(\De(x)),\qquad \forall \
\ga_1,\ga_2\in (U_{r,s}^{\geq0})^*, \  x\in U_{r,s}^{\geq0}.
$$
The identity element is given by $\ep$. We define the linear
functionals $\ga_{\mu},\,\xi_i\in (U_{r,s}^{\geq0})^*$ for any
$i\in I,\ \mu\in Q$ by
\begin{equation*}
\begin{split}
\ga_{\mu}(x\om_{\nu})&=r^{\lan\mu,\nu\ran}s^{-\lan\nu,\mu\ran}\ep(x),\qquad(\forall\  x\in U_{r,s}^+,\ \nu\in Q),\\
\xi_i(x\om_{\nu})&=0,\quad\quad(\forall\  x\in U_{r,s}^{+\be},\ \be\in Q^+- \{\al_i\},\ \nu\in Q) ,\\
\xi_i(e_i\om_{\nu})&=\frac{1}{s_i-r_i},\qquad(\forall\  \nu\in Q).
\end{split}
\end{equation*}
Define a linear map
$$
\phi:U_{r,s}^{\leq0}\lra (U_{r,s}^{\geq0})^*
$$
by $\phi(\om_{\mu}')=\ga_{\mu},\ \phi(f_i)=\xi_i$ and extending it
algebraically. It is straightforward to check that $\phi$ is
well-defined. Now we can define the pairing
$$(\ ,\ ): \ U_{r,s}^{\leq0}\times U_{r,s}^{\geq0}\lra \bQ(r,s)$$
by $(x,y):=\phi(x)(y)$, for any $x\in U_{r,s}^{\leq0},\ y\in
U_{r,s}^{\geq0}$. The condition
$$
(yy',\,x)  =(y\ot y',\,\De(x)),\quad (y, xx')=(\Delta(y), x'\otimes
x),
$$
for $x\in U_{r,s}^{\geq0}, \ y,\,y'\in U_{r,s}^{\leq0}\,$ can be
proved by induction. The remaining conditions are obvious. Moreover,
it is clear that the bilinear form $(\ ,\ )$ is uniquely determined.
\end{proof}

Based on Proposition 2.4, similarly to the proof of Theorem 2.5 in
\cite{BGH1}, we have
\begin{cor}
$U_{r,s}(\mg)$ can be realized as a Drinfel'd double of Hopf
subalgebras $U_{r,s}^{\geq0}$ and $U_{r,s}^{\geq0}$ with respect to
the pairing $(\ ,\ )$, that is,
$$
U_{r,s}(\mg)\cong \D(U_{r,s}^{\geq0},U_{r,s}^{\leq0}).
$$
\end{cor}

As a consequence of the Drinfel'd double structure, with the same
argument of Corollary 2.6 in \cite{BGH1}, we have
\begin{cor}
$U_{r,s}(\mg)$ has the standard triangular decomposition
$$
U_{r,s}(\frak g)\cong U_{r,s}^-\ot U_{r,s}^0\ot U_{r,s}^+,
$$
where $U_{r,s}^0=\bigoplus_{\mu,\nu\in
Q}\mathbb{Q}(r,s)\,\om_{\nu}'\om_{\mu}$ and $
U_{r,s}^{\pm}=\bigoplus_{\be\in Q+}U_{r,s}^{\pm\be}$.
\end{cor}

\medskip
\section{Cocycle deformations of $U_{r,s}^+$}
\medskip

\begin{lem}[\cite{AST}] \ Let $A=\bigoplus_{g\in G}A_g$ be a $G$-graded associative algebra over a field $k$,
where $G$ is an abelian group. Let $\psi: G\times G\to k^* $ be a
$2$-cocycle of the group $G$. We introduce a new multiplication $*$
on $A$ as follows: For any $x\in A_g, \,y\in A_h$, where $g,\,h\in
G$, we define
$$
x*y=\psi(g,h)\,x\,y.
$$
Denote this new algebra by $A^{\psi}$. Then $A^{\psi}$ is a
$G$-graded associative algebra, owing to $\psi$ being a $2$-cocycle.
The algebra  $A^{\psi}$ is called a cocycle deformation of the
algebra $A$ by $\psi$.
\end{lem}

In this section, let us take the parameters $r, s, r', s'\in\bC^*$
and consider both algebras $U_{r,s}^+$ and $U_{r',s'}^+$ to be
defined over the field $\bC$ of complex numbers. Note that both
algebras are $Q$-graded. In view of Lemma 3.1, the argument of the
following main result is interesting.

\begin{theo}
$U_{r,s}^+$ and $U_{r',s'}^+$ are $2$-cocycle deformations of each
other if\,  $rs^{-1}=r's'^{-1}$ or \,$rs^{-1}=r'^{-1}s'$.
\end{theo}
\begin{proof} \
(I) Assume that $rs^{-1}=r's'^{-1}$. In this case, we define a new
product $*$ on $U^+_{r,s}$ as follows
$$
x*y=\psi(\mu,\nu)\,x\,y=(r^{-1}r')^{\lan\mu,\nu\ran}\, x\,y,
$$
for any $x\in U_{r,s}^{+\mu},\ y\in U_{r,s}^{+\nu}$, where $\psi:
Q\times Q\lra\bC^*$ such that
$\psi(\mu,\nu)=(r^{-1}r')^{\lan\mu,\nu\ran}$.

Note that $\psi$ is a bicharacter on $Q\times Q$, which is obviously
a $2$-cocycle of the abelian group $Q$. This fact ensures that the
new $*$-product is associative.

In what follows, it suffices to prove the relations below:
\begin{gather*}
\sum_{k=0}^{1-a_{ij}} (-1)^{k} \binom{1-a_{ij}}{k}_{r_is_i^{-1}}
c^{(k)}_{ij} e_{i}^{*(1-a_{ij}-k)}*e_{j}*e_{i}^{*k} =0, \qquad (i\neq j);\\
 \sum_{k=0}^{1-a_{ij}} (-1)^{k} \binom{1-a_{ij}}{k}_{r_is_i^{-1}}
c^{(k)}_{ij} f_{i}^{*k}*f_{j}*f_{i}^{*(1-a_{ij}-k)} =0, \qquad
(i\neq j),
\end{gather*}
where
$$
c^{(k)}_{ij}=(r_is_i^{-1})^{\frac{k(k-1)}{2}} r^{k\lan
j,i\ran}s^{-k\lan i,j\ran},\qquad (i\neq j).
$$

By the definition of $*$-product, we have
$$
e_{i}^{*(m-k)}*e_{j}* e_{i}^{*k} =(s^{-1}s')^{\frac{m(m-1)}{2}\lan
i,i\ran+(m-k)\lan i,j\ran+k\lan j,i\ran}e_{i}^{m-k} e_{j}
e_{i}^{k}.
$$

Case ($1$): \ $i<j$, \ i.e., $\langle j,i\rangle=0$: when
$m=1-a_{ij}$, we have
\begin{equation*}
\begin{split}
&\sum_{k=0}^{m} (-1)^{k} \binom{m}{k}_{r_is_i^{-1}} c^{(k)}_{ij}
e_{i}^{*(m-k)}* e_{j}*
e_{i}^{*k}\\
&\ =\sum_{k=0}^{m} (-1)^{k} \binom{m}{k}_{r_is_i^{-1}}
c^{(k)}_{ij}(s^{-1}s')^{\frac{m(m-1)}{2}\lan i,i\ran+(m-k)\lan
i,j\ran+k\lan j,i\ran}e_{i}^{m-k}
e_{j} e_{i}^{k}\\
&\ =\sum_{k=0}^{m} (-1)^{k} \binom{m}{k}_{r_is_i^{-1}}
(r_is_i^{-1})^{\frac{k(k-1)}{2}}s^{-k\lan
i,j\ran}(s^{-1}s')^{\frac{m(m-1)}{2}d_i+(m-k)d_ia_{ij}}e_{i}^{m-k}
e_{j} e_{i}^{k}\\
&\ =(s_i^{-1}s_i')^{\frac{m(m-1)}{2}+ma_{ij}}\sum_{k=0}^{m}
(-1)^{k} \binom{m}{k}_{r_is_i^{-1}}
(r_is_i^{-1})^{\frac{k(k-1)}{2}}s_i^{-ka_{ij}}(s_is_i'^{-1})^{ka_{ij}}e_{i}^{m-k}
e_{j} e_{i}^{k}\\
&\ =(s_i^{-1}s_i')^{\frac{m(m-1)}{2}+ma_{ij}}\sum_{k=0}^{m}
(-1)^{k} \binom{m}{k}_{r_is_i^{-1}}
(r_is_i^{-1})^{\frac{k(k-1)}{2}}(s_i'^{-1})^{ka_{ij}}e_{i}^{m-k}
e_{j} e_{i}^{k}\\
&\ =(s_is_i'^{-1})^{\frac{m(m-1)}{2}}\sum_{k=0}^{m} (-1)^{k}
\binom{m}{k}_{r_i's_i'^{-1}}
(r'_is_i'^{-1})^{\frac{k(k-1)}{2}}r'^{k\lan j,i\ran}s'^{-k\lan
i,j\ran} e_{i}^{m-k}
e_{j} e_{i}^{k}\\
&\ =(s_is_i'^{-1})^{\frac{m(m-1)}{2}}\underbrace{\sum_{k=0}^{m}
(-1)^{k} \binom{m}{k}_{r_i's_i'^{-1}} c_{ij}'^{(k)}\, e_{i}^{m-k}
e_{j} e_{i}^{k}}\\
&\ =0.\qquad\qquad\qquad(\textit{$(r',s')$-Serre relations in
$U_{r',s'}^+$})
\end{split}
\end{equation*}

Case ($2$): \ $i>j$, \ i.e., $\langle i,j\rangle=0$: when
$m=1-a_{ij}$, we have
\begin{equation*}
\begin{split}
&\sum_{k=0}^{m} (-1)^{k} \binom{m}{k}_{r_is_i^{-1}} c^{(k)}_{ij}
e_{i}^{*(m-k)}* e_{j}*
e_{i}^{*k}\\
&\ =\sum_{k=0}^{m} (-1)^{k} \binom{m}{k}_{r_is_i^{-1}}
c^{(k)}_{ij}(s^{-1}s')^{\frac{m(m-1)}{2}\lan i,i\ran+(m-k)\lan
i,j\ran+k\lan j,i\ran}e_{i}^{m-k}
e_{j} e_{i}^{k}\\
&\ =\sum_{k=0}^{m} (-1)^{k} \binom{m}{k}_{r_is_i^{-1}}
(r_is_i^{-1})^{\frac{k(k-1)}{2}}r_i^{ka_{ij}}(s^{-1}s')^{\frac{m(m-1)}{2}d_i+kd_ia_{ij}}e_{i}^{m-k}
e_{j} e_{i}^{k}\\
&\ =(s_i^{-1}s_i')^{\frac{m(m-1)}{2}}\sum_{k=0}^{m} (-1)^{k}
\binom{m}{k}_{r_is_i^{-1}}
(r_is_i^{-1})^{\frac{k(k-1)}{2}}r_i^{ka_{ij}}(s_i^{-1}s_i')^{ka_{ij}}e_{i}^{m-k}
e_{j} e_{i}^{k}\\
&\ =(s_i^{-1}s_i')^{\frac{m(m-1)}{2}}\sum_{k=0}^{m} (-1)^{k}
\binom{m}{k}_{r_is_i^{-1}}
(r_is_i^{-1})^{\frac{k(k-1)}{2}}(r_i')^{ka_{ij}}e_{i}^{m-k}
e_{j} e_{i}^{k}\\
&\ =(s_i^{-1}s_i')^{\frac{m(m-1)}{2}}\sum_{k=0}^{m} (-1)^{k}
\binom{m}{k}_{r_i's_i'^{-1}}
(r'_is_i'^{-1})^{\frac{k(k-1)}{2}}r'^{k\lan j,i\ran}s'^{-k\lan
i,j\ran} e_{i}^{m-k}
e_{j} e_{i}^{k}\\
&\ =(s_i^{-1}s_i')^{\frac{m(m-1)}{2}}\sum_{k=0}^{m} (-1)^{k}
\binom{m}{k}_{r_i's_i'^{-1}} c_{ij}'^{(k)} e_{i}^{m-k} e_{j}
e_{i}^{k}=0.
\end{split}
\end{equation*}
Hence, $U_{r,s}^+$ and $U_{r',s'}^+$ are $2$-cocycle deformations of
each other.

(II) Assume that $rs^{-1}=r'^{-1}s'$. In this case, we can define
another new product $*$ on $U^+_{r,s}$ by
$x*y=\psi(\mu,\nu)\,x\,y=(r's^{-1})^{\lan\mu,\nu\ran}\,x\,y$, for
any $x\in U_{r,s}^{+\mu},\ y\in U_{r,s}^{+\nu}$, where $\psi:
Q\times Q\lra\bC^*$ such that
$\psi(\mu,\nu)=(r's^{-1})^{\lan\mu,\nu\ran}$. Thus we have the
following

Case ($1'$): \ $i<j$, \ i.e., $\langle j,i\rangle=0$: when
$m=1-a_{ij}$, we have
\begin{equation*}
\begin{split}
&\sum_{k=0}^{m} (-1)^{k} \binom{m}{k}_{r_is_i^{-1}} c^{(k)}_{ij}
e_{i}^{*(m-k)}* e_{j}*
e_{i}^{*k}\\
&\ =\sum_{k=0}^{m} (-1)^{k} \binom{m}{k}_{r_is_i^{-1}}
c^{(k)}_{ij}(r's^{-1})^{\frac{m(m-1)}{2}\lan i,i\ran+(m-k)\lan
i,j\ran+k\lan j,i\ran}e_{i}^{m-k}
e_{j} e_{i}^{k}\\
&\ =\sum_{k=0}^{m} (-1)^{k} \binom{m}{k}_{r_is_i^{-1}}
(r_is_i^{-1})^{\frac{k(k-1)}{2}}s^{-k\lan
i,j\ran}(r's^{-1})^{\frac{m(m-1)}{2}d_i+(m-k)d_ia_{ij}}e_{i}^{m-k}
e_{j} e_{i}^{k}\\
&\ =(r_i's_i^{-1})^{\frac{m(m-1)}{2}+ma_{ij}}\sum_{k=0}^{m}
(-1)^{k} \binom{m}{k}_{r_is_i^{-1}}
(r_is_i^{-1})^{\frac{k(k-1)}{2}}s_i^{-ka_{ij}}(s_ir_i'^{-1})^{ka_{ij}}e_{i}^{m-k}
e_{j} e_{i}^{k}\\
&\ =(r_i'^{-1}s_i)^{\frac{m(m-1)}{2}}\sum_{k=0}^{m} (-1)^{k}
\binom{m}{k}_{r_i'^{-1}s_i'}
(r_i'^{-1}s_i')^{\frac{k(k-1)}{2}}(r_i'^{-1})^{ka_{ij}}e_{i}^{m-k}
e_{j} e_{i}^{k}\\
&\ =(r_i'^{-1}s_i)^{\frac{m(m-1)}{2}}\sum_{k=0}^{m} (-1)^{k}
\binom{m}{k}_{r_i's_i'^{-1}}
(r_i's_i'^{-1})^{k^2-km}(r_i'^{-1}s_i')^{\frac{k(k-1)}{2}}r_i'^{-ka_{ij}}e_{i}^{m-k}
e_{j} e_{i}^{k}\\
&\ =(r_i'^{-1}s_i)^{\frac{m(m-1)}{2}}\sum_{k=0}^{m} (-1)^{k}
\binom{m}{k}_{r_i's_i'^{-1}}
(r'_is_i'^{-1})^{\frac{k(k-1)}{2}}r'^{k\lan j,i\ran}s'^{-k\lan
i,j\ran} e_{i}^{m-k}
e_{j} e_{i}^{k}\\
&\ =0.
\end{split}
\end{equation*}

Case ($2'$): \ $i>j$, \ i.e., $\langle i,j\rangle=0$: when
$m=1-a_{ij}$, we have
\begin{equation*}
\begin{split}
&\sum_{k=0}^{m}(-1)^{k} \binom{m}{k}_{r_is_i^{-1}} c^{(k)}_{ij}
e_{i}^{*(m-k)}* e_{j}*
e_{i}^{*k}\\
&\ =\sum_{k=0}^{m} (-1)^{k} \binom{m}{k}_{r_is_i^{-1}}
c^{(k)}_{ij}(r's^{-1})^{\frac{m(m-1)}{2}\lan i,i\ran+(m-k)\lan
i,j\ran+k\lan j,i\ran}e_{i}^{m-k}
e_{j} e_{i}^{k}\\
&\ =\sum_{k=0}^{m} (-1)^{k} \binom{m}{k}_{r_is_i^{-1}}
(r_is_i^{-1})^{\frac{k(k-1)}{2}}r_i^{ka_{ij}}(r^{-1}s')^{\frac{m(m-1)}{2}d_i+kd_ia_{ij}}e_{i}^{m-k}
e_{j} e_{i}^{k}\\
&\ =(r_i^{-1}s_i')^{\frac{m(m-1)}{2}}\sum_{k=0}^{m} (-1)^{k}
\binom{m}{k}_{r_is_i^{-1}}
(r_is_i^{-1})^{\frac{k(k-1)}{2}}r_i^{ka_{ij}}(r_i^{-1}s_i')^{ka_{ij}}e_{i}^{m-k}
e_{j} e_{i}^{k}\\
&\ =(r_i^{-1}s_i')^{\frac{m(m-1)}{2}}\sum_{k=0}^{m} (-1)^{k}
\binom{m}{k}_{r_i'^{-1}s_i'}
(r_i'^{-1}s_i')^{\frac{k(k-1)}{2}}(s_i')^{ka_{ij}}e_{i}^{m-k}
e_{j} e_{i}^{k}\\
&\ =(r_i^{-1}s_i')^{\frac{m(m-1)}{2}}\sum_{k=0}^{m} (-1)^{k}
\binom{m}{k}_{r_i's_i'^{-1}}(r_i's_i'^{-1})^{k^2-km}
(r_i'^{-1}s_i')^{\frac{k(k-1)}{2}}(s_i')^{ka_{ij}}e_{i}^{m-k}
e_{j} e_{i}^{k}\\
&\ =(r_i^{-1}s_i')^{\frac{m(m-1)}{2}}\sum_{k=0}^{m} (-1)^{k}
\binom{m}{k}_{r_i's_i'^{-1}}
(r'_is_i'^{-1})^{\frac{k(k-1)}{2}}r'^{k\lan j,i\ran}s'^{-k\lan
i,j\ran} e_{i}^{m-k}
e_{j} e_{i}^{k}\\
&\ =0.
\end{split}
\end{equation*}
Hence, $U_{r,s}^+$ and $U_{r',s'}^+$ are $2$-cocycle deformations of
each other.
\end{proof}

\begin{cor}
$(\text{\rm i})$ \ $U_{r,s}^+ $ and $U_{s^{-1},r^{-1}}^+$ $($the
so-called associated object of the former in \cite{BGH1}$)$ are
$2$-cocycle deformations of each other. Moreover, $U_{r,s}^+
=U_{s^{-1},r^{-1}}^+$ if and only if $\,rs=1$.

$(\text{\rm ii})$ \ In particular, if $\,rs^{-1}=q^2$, $U_{r,s}^+$,
$U_{q^2,1}^+$ and $U_{q,q^{-1}}^+$ are $2$-cocycle deformations of
each other.\hfill\qed
\end{cor}

\medskip
\section{Realization and Kashiwara's skew differential operators}
\medskip

The following result arises from Kashiwara's work \cite{Ka} (in
one-parameter case).
\begin{pro}
For $P \in U_{r,s}^+$, there exist unique $L,\,R \in U_{r,s}^+$
satisfying the following equation
$$
[\,P,f_i\,]=\frac{\om_i L - \om_i'R}{r_i-s_i},
$$
where we define $ \pa_i(P)= L$ and $\pa_i'(P)= R$.
\end{pro}
\begin{proof} \ Assume that
$$
\frac{\om_i L_1 - \om_i'R_1}{r_i-s_i}=\frac{\om_i L_2 -
\om_i'R_2}{r_i-s_i},
$$
then we have
$$
\om_i(L_1-L_2)-\om_i'(R_1-R_2)=0.
$$
Using the triangular decomposition of $U_{r,s}(\mg)$ in Corollary
2.6, we get $L_1=L_2,\, R_1=R_2$. This means that the uniqueness of
$L$ and $R$ is clear.

To show the existence of $L$ and $R$, we consider its graded
decomposition $U_{r,s}^{+}=\bigoplus_{\nu\in Q^{+}}U_{r,s}^{+\nu}$.
We proceed to prove this by induction on the height of weights,
namely, $\text{ht}(\nu)=\sum m_i$, for $\nu=\sum m_i\al_i\;(m_i\in
\mathbb{Z}_{\geq 0})$. If $\text{ht}{(\nu)}=1$ and $P=e_j$, then we
choose $L=R=\delta_{i,j}$. Suppose that for $P\in U_{r,s}^{+\nu}$,
there exist $L,\, R$ satisfying
$$
[\,P,f_i\,]=\frac{\om_i L - \om_i'R}{r_i-s_i}.
$$
Then for $e_j\,P\in U_{r,s}^{+(\nu+1)}$ with
$\textrm{ht}(e_j\,P)=\nu+1$, we have
\begin{equation*}
\begin{split}
[\,e_j&\,P, f_i\,]=e_j\,[\,P,f_i\,]+[\,e_j,f_i\,]\,P\\
&=\frac{1}{r_i-s_i}\Bigl\{\om_i\Bigl(r^{-\lan j, i\ran}s^{\lan i,
j\ran}e_j L+\delta_{i,j}P\Bigr)-\om_i'\Bigl(r^{\lan
i,j\ran}s^{-\lan j,i\ran} e_j
R+\delta_{i,j}P\Bigr)\Bigr\}\\
&=\frac{1}{r_i-s_i}\Bigl\{\om_i\Bigl(r^{-\lan j, i\ran}s^{\lan i,
j\ran}e_j \pa_i(P)+\delta_{i,j}P\Bigr)-\om_i'\Bigl(r^{\lan
i,j\ran}s^{-\lan j,i\ran} e_j \pa'_i(P)+\delta_{i,j}P\Bigr)\Bigr\},
\end{split}
\end{equation*}
by the induction hypothesis. In particular, we get
\begin{equation*}
\begin{split}
\pa_i(e_jP)&=r^{-\lan j, i\ran}s^{\lan i, j\ran}e_j
\pa_i(P)+\delta_{i,j}P,\\
\pa_i'(e_jP)&=r^{\lan i,j\ran}s^{-\lan j,i\ran} e_j
\pa_i'(P)+\delta_{i,j}P.
\end{split}\tag{*}
\end{equation*}

This completes the proof.
\end{proof}

From (*), we easily get
$$
\partial_i(e_i^m)=(m)_{r_i^{-1}s_i}e_i^{m-1},\qquad
\partial_i'(e_i^m)=(m)_{r_is_i^{-1}}e_i^{m-1}.$$

For $i\in I$, we introduce the operator $E_i:\
U_{r,s}^+\longrightarrow U_{r,s}^+$ defined by
$$E_i\,x=e_i\,x,\quad \text{\it for any} \ \, x\in U_{r,s}^+.$$
Then the following lemma can be proved inductively by using
Proposition 4.1.
\begin{lem}
For $m\in\bZ_{+}$, $i,\, j\in I$, the following commutation
relations hold
\begin{gather*}
\pa_i\,\pa_j'=r^{\lan j,i\ran}s^{-\lan i,j\ran}\pa_j'\,\pa_i,\\
\pa_i^m\,E_j=r^{-m\lan j, i\ran}s^{m\lan i, j\ran}\,E_j\,\pa_i^m
+\de_{i,j}(m)_{r_i^{-1}s_i}\pa_i^{m-1},\\
\pa_i'^m\,E_j=r^{m\lan i,j\ran}s^{-m\lan j,i\ran}\,E_j\,\pa_i'^m
+\de_{i,j}(m)_{r_is_i^{-1}}\pa_i'^{m-1}.
\end{gather*}
\end{lem}

\begin{theo} For any $i\neq j\in I$, we have
\begin{gather*}
\sum_{k=0}^{1-a_{ij}} (-1)^{k} \binom{1-a_{ij}}{k}_{r_is_i^{-1}}
c^{(k)}_{ij}\, \pa_i^{k}\,\pa_j\,
\pa_i^{1-a_{ij}-k} =0,\tag{\textrm{i}}\\
\sum_{k=0}^{1-a_{ij}} (-1)^{k} \binom{1-a_{ij}}{k}_{r_is_i^{-1}}
c^{(k)}_{ij}\, \pa_i'^{k}\,\pa_j'\, \pa_i'^{1-a_{ij}-k}
=0,\tag{\textrm{ii}}
\end{gather*}
which give rise to two $($operators\,$)$ realizations of $U^+_{r,s}$
via assigning the generators $e_i$'s of $U_{r,s}^{+}$ to the
Kashiwara's skew differential operators $\partial_i$'s or
 $\partial_i'$'s respectively.
\end{theo}
\begin{proof} \ (i) \ For any $u\in U_{r,s}^{+\mu}$, we will prove
the formula
$$
\sum_{k=0}^{1-a_{ij}} (-1)^{k} \binom{1-a_{ij}}{k}_{r_is_i^{-1}}
c^{(k)}_{ij}\,
\pa_i^{1-a_{ij}-k}\,\pa_j\,\pa_i^k\,u=0,\leqno(\text{$**$})
$$
by induction on $\text{ht}(u)=\mu$.  For any $v\in U_{r,s}^{+\nu}$
with $\text{ht}(v)=\nu< \text{ht}(u)=\mu$, we assume that $(**)$
holds. Write $u$ as $u=e_\ell\, v=E_\ell\, v$, for some $\ell\in I$
and some $v\in U_{r,s}^{+\nu}$. Now put $m=1-a_{ij}$.

First, we note that
\begin{equation*}
\begin{split}
&\pa_i^{m-k}\, \pa_j\,\pa_i^{k}\,E_\ell \\
&\quad=\pa_i^{m-k}\, \pa_j
\left\{\,r^{-k\lan \ell , i\ran}s^{k\lan i, \ell \ran}\,E_\ell \,\pa_i^k +\de_{i,\ell }\,(k)_{r_i^{-1}s_i}\,\pa_i^{k-1}\,\right\}\\
&\quad= \pa_i^{m-k} \left\{\,r^{-k\lan \ell , i\ran}s^{k\lan i,
\ell \ran}\Bigl(r^{-\lan \ell , j\ran}s^{\lan j, \ell
\ran}\,E_\ell \,\pa_j+\de_{\ell ,j}\Bigr)\,\pa_i^k
+\de_{i,\ell }\,(k)_{r_i^{-1}s_i}\,\pa_j\,\pa_i^{k-1}\,\right\}\\
&\quad=
r^{-k\lan \ell , i\ran}s^{k\lan i, \ell \ran}r^{-\lan \ell , j\ran}s^{\lan j, \ell \ran}\times\\
&\qquad\left\{r^{(k-m)\lan \ell , i\ran}
s^{(m-k)\lan i, \ell \ran}\,E_\ell \,\pa_i^{m-k}+\de_{\ell ,i}\,(m-k)_{r_i^{-1}s_i}\,\pa_i^{m-k-1}\right\}\pa_j\,\pa_i^k\\
&\qquad+\de_{\ell ,j}\,r^{-k\lan \ell , i\ran}s^{k\lan i, \ell
\ran}\,{\pa_i}^{m-k}\,\pa_i^k
+\de_{i,\ell }\,(k)_{r_i^{-1}s_i}\,{\pa_i}^{m-k}\,\pa_j\,\pa_i^{k-1}\\
&\quad=
r^{-m\lan \ell , i\ran}s^{m\lan i, \ell \ran}\,r^{-\lan \ell , j\ran}s^{\lan j,\ell \ran}\,E_\ell \,\pa_i^{m-k}\,\pa_j\,\pa_i^k\\
&\qquad+ \de_{\ell ,i}\,r^{-k\lan \ell , i\ran}s^{k\lan i, \ell \ran}r^{-\lan \ell , j\ran}s^{\lan j, \ell \ran}\,(m-k)_{r_i^{-1}s_i}\,\pa_i^{m-k-1}\,\pa_j\,\pa_i^k\\
&\qquad+\de_{i,\ell
}\,(k)_{r_i^{-1}s_i}\,{\pa_i}^{m-k}\,\pa_j\,\pa_i^{k-1}+\de_{\ell
,j}\,r^{-k\lan \ell , i\ran}s^{k\lan i, \ell \ran}\,{\pa_i}^{m}.
\end{split}
\end{equation*}
Consequently, we obtain
\begin{equation*}
\begin{split}
&\sum_{k=0}^{m} (-1)^{k} \binom{m}{k}_{r_is_i^{-1}}\,
c^{(k)}_{ij}\, \pa_i^{k}\,\pa_j\,
\pa_i^{m-k}\,E_\ell \\
&\ =  r^{-m\lan \ell , i\ran}s^{m\lan i, \ell \ran}r^{-\lan \ell ,
j\ran}s^{\lan j,\ell \ran}\,E_\ell \,\sum_{k=0}^{m} (-1)^{k}
\binom{m}{k}_{r_is_i^{-1}}\,
c^{(k)}_{ij}\,\pa_i^{m-k}\,\pa_j\,\pa_i^k\\
&\quad+\de_{\ell ,i}\,r^{-\lan i, j\ran}s^{\lan j,
i\ran}\sum_{k=0}^{m} (-1)^{k}
\binom{m}{k}_{r_is_i^{-1}}(m-k)_{r_i^{-1}s_i}
c^{(k)}_{ij}r^{-k\lan i, i\ran}s^{k\lan i, i\ran}\pa_i^{m-k-1}\,\pa_j\,\pa_i^k\\
&\quad+\de_{i,\ell }\,\sum_{k=0}^m (-1)^{k}
\binom{m}{k}_{r_is_i^{-1}}\,(k)_{r_i^{-1}s_i}\,
c^{(k)}_{ij}\,{\pa_i}^{m-k}\,\pa_j\,\pa_i^{k-1}\\
&\quad+\de_{\ell
,j}\,{\pa_i}^{m}\,\sum_{k=0}^{m} (-1)^{k}
\binom{m}{k}_{r_is_i^{-1}}\,
c^{(k)}_{ij}\,r^{-k\lan j, i\ran}s^{k\lan i, j\ran}\\
&\ =S_1+S_2+S_3+S_4=S_1,
\end{split}
\end{equation*}
where $S_2=-S_3$, $S_4=0$ (by Lemma 4.4 below), and
\begin{gather*}
S_1= r^{-m\lan \ell , i\ran}s^{m\lan i, \ell \ran}r^{-\lan \ell ,
j\ran}s^{\lan j,\ell \ran}\,E_\ell \,\sum_{k=0}^{m} (-1)^{k}
\binom{m}{k}_{r_is_i^{-1}}\,
c^{(k)}_{ij}\,\pa_i^{m-k}\,\pa_j\,\pa_i^k,\\
S_2=\de_{\ell ,i}\sum_{k=0}^{m-1} (-1)^{k}
\binom{m}{k}_{r_is_i^{-1}}(m-k)_{r_i^{-1}s_i}
c^{(k)}_{ij}r^{-\lan i, j\ran}s^{\lan j, i\ran}r^{-k\lan i, i\ran}s^{k\lan i, i\ran}\pa_i^{m-k-1}\pa_j\pa_i^k,\\
S_3=\de_{i,\ell }\,\sum_{k=1}^{m} (-1)^{k}
\binom{m}{k}_{r_is_i^{-1}}\,(k)_{r_i^{-1}s_i}\,
c^{(k)}_{ij}\,{\pa_i}^{m-k}\,\pa_j\,\pa_i^{k-1},
\\
S_4=\de_{\ell ,j}\,{\pa_i}^{m}\,\sum_{k=0}^{m} (-1)^{k}
\binom{m}{k}_{r_is_i^{-1}}\, c^{(k)}_{ij}\,r^{-k\lan j,
i\ran}s^{k\lan i, j\ran}.
\end{gather*}

Now according to the inductive hypothesis, we get $S_1\, v=0$. So we
proved the equality $(**)$. This means that
 $U_{r,s}^{+}$ can be realized by identifying the generators $e_i$'s
 with the skew differential
operators $\pa_i$'s, that is, the algebra generated by the
$\partial_i$'s is a homomorphic image of $U_{r,s}^{+}$.

(ii) \ Similarly, we can prove another identity (ii), which shows
that $U_{r,s}^{+}$ can be realized by the skew differential
operators $\pa_i'$'s.
\end{proof}

\begin{lem} \ $(\textrm{\rm i})$ \ $S_2=-S_3$.

$(\textrm{\rm ii})$ \ $S_4=0$.
\end{lem}
\begin{proof}\ (ii) follows from the identity below:
$$
\sum_{k=0}^{m} (-1)^{k} \binom{m}{k}_{r_is_i^{-1}}
c^{(k)}_{ij}r^{-k\lan j, i\ran}s^{k\lan i, j\ran}=\sum_{k=0}^{m}
(-1)^{k}
\binom{m}{k}_{r_is_i^{-1}}(r_is_i^{-1})^{\frac{k(k-1)}{2}}=0.
$$

As for (i), we notice that $(n)_{q^{-1}}=q^{1-n}(n)_q$,
$1-m=a_{ij}$, and
\begin{equation*}
\begin{split}
c_{ij}^{(k-1)}&=c_{ij}^{(k)}r^{-\lan
j,i\ran}s^{\lan i,j\ran}(r_is_i^{-1})^{1-k},\\
r^{-\lan j,i\ran-\lan i,j\ran}s^{\lan i,j\ran+\lan
j,i\ran}&=(r_i^{-1}s_i)^{a_{ij}}=(r_is_i^{-1})^{m-1},\\
 \binom{m}{k{-}1}_{r_is_i^{-1}}
(m{-}k{+}1)_{r_i^{-1}s_i}&=\binom{m}{k{-}1}_{r_is_i^{-1}}
(m{-}k{+}1)_{r_is_i^{-1}}(r_is_i^{-1})^{k-m}\\
&=\binom{m}{k}_{r_is_i^{-1}}
(k)_{r_is_i^{-1}}(r_is_i^{-1})^{k-m}\\
&=\binom{m}{k}_{r_is_i^{-1}}
(k)_{r_i^{-1}s_i}(r_is_i^{-1})^{2k-m-1},
\end{split}
\end{equation*}
so that
$$
c_{ij}^{(k-1)}\,r^{-\lan i, j\ran}s^{\lan j,
i\ran}(r_i^{-1}s_i)^{k-1}=c_{ij}^{(k)}(r_is_i^{-1})^{m-1+2(1-k)}.$$
Consequently, we obtain
\begin{equation*}
\begin{split}
S_2&=\de_{\ell ,i}\sum_{k=0}^{m-1} (-1)^{k}
\binom{m}{k}_{r_is_i^{-1}}(m{-}k)_{r_i^{-1}s_i}
c^{(k)}_{ij}\,r^{-\lan i, j\ran}s^{\lan j,
i\ran}(r_i^{-1}s_i)^k\,\pa_i^{m{-}k{-}1}\,\pa_j\,\pa_i^k\\
&=\de_{\ell ,i}\sum_{k=1}^{m} (-1)^{k{+}1}
\binom{m}{k{-}1}_{r_is_i^{-1}}(m{-}k{+}1)_{r_i^{-1}s_i}
c^{(k-1)}_{ij}\,r^{-\lan i, j\ran}s^{\lan j,
i\ran}(r_i^{-1}s_i)^{k-1}\times \\
&\hskip9cm\times\pa_i^{m-k}\,\pa_j\,\pa_i^{k-1}\\
&=-\de_{\ell ,i}\sum_{k=1}^{m}
(-1)^k\binom{m}{k}_{r_is_i^{-1}}\,(k)_{r_i^{-1}s_i}\,
c^{(k)}_{ij}\,{\pa_i}^{m-k}\,\pa_j\,\pa_i^{k-1}\\
&=-S_3.
\end{split}
\end{equation*}

This completes the proof.
\end{proof}

\vskip30pt \centerline{\bf ACKNOWLEDGMENT}

\vskip15pt N.H. thanks M. Rosso and V.K. Dobrev for their useful
comments on \cite{BGH1} \& \cite{BGH2}, as well as the ICTP for its
hospitality and support when he visited the ICTP Mathematics Group
from March to August of 2006, Trieste, Italy.

\bigskip
\bibliographystyle{amsalpha}

\begin{thebibliography}{9999}
\medskip


\bibitem{AST} Artin M, Schelter W, and Tate J. \textit{Quantum deformations of
$GL(n)$}, Comm. Pure Appl. Math. {\textbf44} (1991), 879--895

\bibitem{Res} Reshetikhin N. \textit{Multiparameter quantum groups and
twisted quasitriangular Hopf algebras}, Lett. Math. Phys. {\textbf
20}, (1990), pp. 331--335

\bibitem{HW} Hu N, Wang X. \textit{Quantizations of the generalized-Witt
algebra and of Jacobson-Witt algebra in the modular case},
arXiv.Math:QA/0602281, J. Algebra,  \textbf{312} (2007), 902--929


\bibitem{BW0} Benkart G and Witherspoon S. \textit{A Hopf structure for
down-up algebras}, Math. Z. \textbf{238} (3), (2001), 523--553

\bibitem{BW1}
Benkart G and Witherspoon S. \textit{Two-parameter quantum groups
and Drinfel'd doubles}, Algebr. Represent. Theory, \textbf{7}
(2004), 261--286


\bibitem{BW2} Benkart G and Witherspoon S. \textit{Representations of two-parameter quantum
groups and Schur-Weyl duality}, Hopf algebras, pp. 65--92, Lecture
Notes in Pure and Appl. Math., \textbf{237}, Dekker, New York, 2004

\bibitem{BW3} Benkart G and Witherspoon S. \textit{Restricted two-parameter quantum groups},
Fields Institute Communications, ``Representations of Finite
Dimensional Algebras and Related Topics in Lie Theory and Geometry",
vol. \textbf{40}, Amer. Math. Soc., Providence, RI, 2004, pp.
293--318

\bibitem{BKL} Benkart G, Kang S-J, Lee K H. \textit{On the center of two-parameter quantum groups},
Proc. Roy. Soc. Edingburg Sect. A, \textbf{136} (3), (2006),
445--472

\bibitem{BGH1} Bergeron N, Gao Y, Hu N. \textit{Drinfel'd doubles and
Lusztig's symmetries of two-parameter quantum groups},
arXiv.Math.RT/0505614, J. of Algebra, \textbf{301} (2006), 378--405

\bibitem{BGH2} Bergeron N, Gao Y, Hu N. \textit{Representations of
two-parameter quantum orthogonal groups and symplectic groups},
arXiv.Math.QA/0510124, AMS/IP Studies in Advanced Mathematics,
\textbf{39} (2007), 1--21

\bibitem{HS} Hu N, Shi Q. \textit{The two-parameter quantum group of exceptional type $G_2$ and
Lusztig's symmetries}, arXiv:Math.QA/0601444, Pacific J. Math. Vol.
230 (2), (2007), 327--346

\bibitem{BH} Bai X, Hu N. \textit{Two-parameter quantum groups of exceptional type E-series
and convex PBW type basis}, arXiv.Math.QA/0605179,  Algebra
Colloquium (in press)

\bibitem{HW2} Hu N, Wang X. \textit{Two-parameter Lusztig's small quantum
groups of type $B$ and their ribbon elements}, Preprint 2006, (42
pages)

\bibitem{CHW} Chen R, Hu N, Wang, X. \textit{Two-parameter Lusztig's small quantum
groups of type $C$ and their ribbon elements}, Preprint 2007, (43
pages)

\bibitem{BH2} Bai X, Hu N. \textit{Two-parameter Lusztig's small quantum
groups of type $D$ and their ribbon elements}, Preprint 2006, (36
pages)

\bibitem{HRZ} Hu N, Rosso M, Zhang H. \textit{Two-parameter affine quantum group $U_{r,s}(\widehat
{\frak{sl}_n})$, Drinfeld realization and quantum affine Lyndon
basis}, Comm. Math. Phys. (in press)

\bibitem{HZ} Hu N, Zhang H. \textit{Vertex representations of two-parameter
quantum affine algebras $U_{r,s}(\widehat{\frak{g}})$: the simply
laced cases}, Preprint 2006, (40 pages)


\bibitem{Ka} Kashiwara M. \textit{On crystal
bases of the q-analogue of universal enveloping algebras}, Duke
Math. J. \textbf{63} (1991), 465--516

\bibitem{Ja} Jantzen J C. \textit{Lectures on Quantum Groups},
Graduate Studies in Mathematics {\bf6}, Amer. Math. Soc. Providence,
1996


\bibitem{Ri} Ringel C. \textit{Hall algebras and quantum groups},
Invent. Math. \textbf{101}  (1990), 583--591



\end{thebibliography}

\end{document}